\documentclass{tac}
\keywords{bicategory, finite products, discrete, comonad, 
Eilenberg-Moore object}
\usepackage{epsfig,epic,eepic,latexsym, amssymb,amsmath}
\usepackage[all]{xy}
\usepackage{enumerate}
\amsclass{18A25}
\copyrightyear{2009}

\title{Bicategories of Spans as Cartesian Bicategories}
\author{Stephen Lack, R.F.C Walters, and R.J. Wood}
\thanks{The authors gratefully acknowledge financial support from 
the Australian ARC, the Italian CNR and the Canadian NSERC.
Diagrams typeset using M. Barr's diagram package, diagxy.tex.}
\address{\\[3pt]
School of Computing and Mathematics\\
University of Western Sydney\\
Locked Bag 1797 South Penrith DC NSW 1797\\
Australia\\
and\\[3pt]
\\Dipartimento di Scienze delle Cultura\\
Politiche e dell'Informazione\\
Universit\`a dell Insubria, Italy\\[3pt] 
and\\[3pt]
Department of Mathematics and Statistics\\
Dalhousie University\\
Halifax, NS, B3H 3J5 Canada}

\eaddress{s.lack@uws.edu.au, robert.walters@uninsubria.it, rjwood@dal.ca}

\newtheorem{axiom}{Axiom}

\let\thm\theorem
\let\prp\proposition
\let\cor\corollary
\let\lem\lemma
\let\dfn\definition
\let\eth\endtheorem
\let\prf\proof
\let\frp\endproof

\let\rmk\remark
\let\axm\axiom

\input diagxy
\xyoption{curve}

\newcommand{\s}{\scalefactor{0.5}}
\newcommand{\qqq}{\scalefactor{0.75}}

\newcommand{\f}{{\kern -.25em}:{\kern -.25em}}

\newcommand{\ra}{{\s\to}}
\newcommand{\tra}{{\s\two}}
\newcommand{\la}{\,{\s\toleft}\,}

\newcommand{\arr}{\mathbf{2}}

\newcommand{\iso}{\cong}
\newcommand{\laj}{\dashv}

\newcommand{\x}{\times}
\newcommand{\ox}{\otimes}

\newcommand{\op}{{^\mathrm{op}}}

\newcommand{\inv}{^{-1}}

\newcommand{\CAT}{\mathbf{CAT}}

\newcommand{\one}{\mathbf{1}}

\newcommand{\dis}{\mathrm{Dis}}

\newcommand{\map}{\mathrm{Map}}

\newcommand{\rel}{\mathrm{Rel}\,}
\newcommand{\spn}{\mathrm{Span}\,}
\newcommand{\com}{\mathrm{Com}}

\newcommand{\E}{{\cal E}}

\newcommand{\bB}{\mathbf{B}}
\newcommand{\bC}{\mathbf{C}}
\newcommand{\bD}{\mathbf{D}}

\newcommand{\bG}{\mathbf{G}}

\newcommand{\bM}{\mathbf{M}}

\begin{document}

\maketitle

\begin{abstract}
Bicategories of spans are characterized as cartesian bicategories in 
which every comonad has an Eilenberg-Moore object and every left 
adjoint arrow is comonadic. 
\end{abstract}

\section{Introduction}\label{intro}
Let $\E$ be a category with finite limits. For the bicategory $\spn\E$,
the locally full subbicategory $\map\spn\E$ determined by the left adjoint
arrows is essentially locally discrete, meaning that each hom category 
$\map\spn\E(X,A)$ is an equivalence relation, and so is equivalent to a 
discrete category. Indeed, a span 
$x\f X\la S\ra A\f a$ has a right adjoint if and only if $x\f S\ra X$ is
invertible. The functors 
$$\map\spn\E(X,A)\to\E(X,A)\quad\mbox{given by}\quad(x,a)\mapsto ax\inv$$
provide equivalences of categories which are the effects on homs for a 
biequivalence 
$$\map\spn\E\ra\E\, .$$ 

Since $\E$ has finite products, $\map\spn\E$ has finite products 
{\em as a bicategory}. We refer the reader to \cite{ckww} for a 
thorough treatment of bicategories with finite products.
Each hom category $\spn\E(X,A)$ is {\em isomorphic} to the slice
category $\E/(X\x A)$ which has binary products given by pullback in
$\E$ and terminal object $1\f X\x A\ra X\x A$. Thus $\spn\E$ is a
{\em precartesian bicategory} in the sense of \cite{ckww}. The
canonical lax monoidal structure 
$$\spn\E\x \spn\E\to\spn\E\toleft\one$$
for this precartesian bicategory is seen to have its binary aspect 
given on arrows by
$$(X\xleftarrow{x} S \xrightarrow{y} A\,,\, Y\xleftarrow{y} T\xrightarrow{b} B )
\mapsto (X\x Y \xleftarrow{x\x y} S\x T \xrightarrow{a\x b} A\x B)\, ,$$
and its nullary aspect provided by
$$1\xleftarrow1 1 \xrightarrow1 1\, ,$$
the terminal object of $\spn\E(1,1)$. Both of these lax functors are
readily seen to be pseudofunctors so that $\spn\E$ is a {\em cartesian
bicategory} as in \cite{ckww}.

The purpose of this paper is to characterize those cartesian bicategories
$\bB$ which are biequivalent to $\spn\E$, for some category $\E$ with 
finite limits. Certain aspects of a solution to the problem are immediate.
A biequivalence $\bB\sim\spn\E$ provides $$\map\bB\sim\map\spn\E\sim\E$$
so that we must ensure firstly that $\map\bB$ is essentially locally discrete.
From the characterization of bicategories of relations as locally ordered
cartesian bicategories in \cite{caw} one suspects that the following axiom
will figure prominently in providing essential local discreteness for
$\map\bB$. 
\axm\label{frob}{\em Frobenius:}\quad A cartesian bicategory $\bB$
is said to satisfy the {\em Frobenius} axiom if, for each $A$ in $\bB$,
$A$ is Frobenius.
\eth
\noindent
Indeed Frobenius objects in cartesian bicategories were defined 
and studied in \cite{ww} where amongst other things it is shown that if $A$ 
is Frobenius in cartesian $\bB$ then, for all $X$, $\map\bB(X,A)$ is a 
groupoid. (This theorem was generalized considerably in \cite{lsw} which 
explained further aspects of the Frobenius concept.) However, essential local
discreteness for $\map\bB$ requires also that the $\map\bB(X,A)$ be 
ordered sets (which is automatic for locally ordered $\bB$). Here we
study also {\em separable} objects in cartesian bicategories for which we
are able to show that if $A$ is  separable in cartesian $\bB$ then, 
for all $X$, $\map\bB(X,A)$ is an ordered set and a candidate axiom is:
\axm\label{sepax}{\em Separable:}\quad A cartesian bicategory $\bB$
is said to satisfy the {\em Separable} axiom if, for each $A$ in $\bB$,
$A$ is separable.
\eth

In addition to essential local discreteness, it is clear that we will
need an axiom which provides {\em tabulation} of each arrow of $\bB$
by a span of maps. Since existence of Eilenberg-Moore objects is a 
basic 2-dimensional limit concept, we will express tabulation in terms of
this requirement; we note that existence of pullbacks in $\map\bB$ 
follows easily from tabulation. In the bicategory $\spn\E$, the comonads
$G\f A\ra A$ are precisely the symmetric spans $g\f A\la X\ra A\f g$;
the map $g\f X\ra A$ together with $g\eta_g\f g\ra gg^*g$ provides an
Eilenberg-Moore coalgebra for $g\f A\la X\ra A\f g$. 
We will posit:
\axm\label{emc}{\em Eilenberg-Moore for Comonads:}\quad Each comonad 
$(A,G)$ in $\bB$ has an Eilenberg-Moore object.
\eth
Conversely, any map (left adjoint) $g\f X\ra A$ in $\spn\E$ provides
an Eilenberg-Moore object for the comonad $gg^*$.
\noindent
We further posit:
\axm\label{mc}{\em Maps are Comonadic:}\quad Each left adjoint $g\f X\ra A$ 
in $\bB$ is comonadic.
\eth
\noindent
from which, in our context, we can also deduce the Frobenius and
Separable axioms. 

In fact we shall also give, in Proposition~\ref{easy} below, a straightforward proof 
that $\map\bB$ is locally essentially discrete whenever Axiom~\ref{mc} holds. 
But due to the importance of the Frobenius and separability conditions in other
contexts, we have chosen to analyze them in their own right.

\section{Preliminaries}\label{prelim}
We recall from \cite{ckww} that a bicategory $\bB$ (always, for
convenience, assumed to be normal) is said to be 
{\em cartesian} if the subbicategory of maps (by which we mean
left adjoint arrows), $\bM=\map\bB$, has finite products $-\x-$ 
and $1$; each hom-category $\bB(B,C)$ has finite products $-\wedge-$
and $\top$; and a certain derived tensor product $-\ox-$ and $I$ 
on $\bB$, extending the product structure of $\bM$, is functorial. 
As in \cite{ckww}, 
we write $p$ and $r$ for the first and second projections at the global level, 
and similarly $\pi$ and $\rho$ for the projections at the local level.
If $f$ is a map of $\bB$ --- an arrow of $\bM$ --- we will write
$\eta_f,\epsilon_f\f f\laj f^*$ for a chosen adjunction in $\bB$
that makes it so. It was shown that the derived tensor product of a 
cartesian bicategory underlies a symmetric monoidal bicategory 
structure. We recall too that in \cite{ww} Frobenius objects in a 
general cartesian bicategory were defined and studied. We will need 
the central results of that paper too. Throughout this paper, $\bB$ 
is assumed to be a cartesian bicategory.

As in \cite{ckww} we write 
$$\bfig
\Atriangle/->`->`/[\bG={\rm Gro}\bB`\bM`\bM;\partial_0`\partial_1`]
\efig$$
for the Grothendieck span corresponding to
$$\bM\op\x\bM\to^{i\op\x i}\bB\op\x\bB\to^{\bB(-,-)}\CAT$$
where $i\f\bM\ra\bB$ is the inclusion. A typical arrow of $\bG$,
$(f,\alpha,u)\f(X,R,A)\ra(Y,S,B)$ can be depicted by a square
\begin{equation}\label{square}
\bfig
\square(0,0)[X`Y`A`B;f`R`S`u]
\morphism(125,250)|m|<250,0>[`;\alpha]
\efig
\end{equation}
and such arrows are composed by pasting. A 2-cell 
$(\phi,\psi)\f(f,\alpha,u)\ra(g,\beta,v)$ in $\bG$ is a pair of
2-cells $\phi\f f\ra g$, $\psi\f u\ra v$ in $\bM$ which satisfy
the obvious equation. The (strict) pseudofunctors $\partial_0$ 
and $\partial_1$ should be regarded as {\em domain} and {\em codomain}
respectively. Thus, applied to (\ref{square}), $\partial_0$ gives $f$ 
and $\partial_1$ gives $u$. The bicategory $\bG$ also has finite products,
which are given on objects by $-\ox-$ and $I$; these are preserved by
$\partial_0$ and $\partial_1$.

The Grothendieck span  can also be thought of as giving a double category 
(of a suitably weak flavour), although we shall not emphasize that point of view.

\subsection{}\label{xredux}
The arrows of $\bG$ are particularly well suited to relating the
various product structures in a cartesian bicategory. In 3.31 of
\cite{ckww} it was shown that the local binary product, for
$R,S\f X{\s \two}A$, can be recovered to within isomorphism from 
the defined tensor product by
$$R\wedge S\cong d^*_A(R\ox S)d_X$$
A slightly more precise version of this is that
the mate of the isomorphism above, with respect to the single
adjunction $d_A\laj d^*_A$, defines an arrow in $\bG$
$$\bfig
\square(0,0)[X`X\ox X`A`A\ox A;d_X`R\wedge S`R\ox S`d_A]
\morphism(125,250)|m|<250,0>[`;]
\efig$$
which when composed with the projections of $\bG$, recovers the
local projections as in
$$
\bfig
\square(0,0)|almb|[X`X\ox X`A`A\ox A;d_X`R\wedge S`R\ox S`d_A]
\morphism(125,250)|m|<250,0>[`;]
\square(500,0)|amrb|[X\ox X`X`A\ox A`A;p_{X,X}`R\ox S`R`p_{A,A}]
\morphism(625,250)|m|<250,0>[`;\tilde p_{R,S}]
\place(1250,250)[\iso]
\square(1550,0)|almb|[X`X`A`A;1_X`R\wedge S`R`1_A]
\morphism(1675,250)|m|<250,0>[`;\pi]
\efig$$
for the first projection, and similarly for the second.
The unspecified $\iso$ in $\bG$ is given by a pair of
convenient isomorphisms $p_{X,X}d_X\iso 1_X$ and $p_{A,A}d_A\iso 1_A$
in $\bM$. Similarly, when $R\wedge S\ra R\ox S$ is composed
with $(r_{X,X},\tilde r_{R,S},r_{A,A})$ the result is 
$(1_X,\rho,1_A)\f R\wedge S\ra S$.

\subsection{}\label{bc}
Quite generally, an arrow of $\bG$ as given by the square (\ref{square})
will be called a {\em commutative} square if $\alpha$ is invertible.
An arrow of $\bG$ will be said to satisfy the 
{\em Beck condition} if the mate
of $\alpha$ under the adjunctions $f\laj f^*$ and $u\laj u^*$, as given
in the square below (no longer an arrow of $\bG$), is invertible.
$$\bfig
\square(1000,0)/<-`->`->`<-/[X`Y`A`B;f^*`R`S`u^*]
\morphism(1125,250)|m|<250,0>[`;\alpha^*]
\efig$$
Thus Proposition 4.7 of \cite{ckww} says that projection squares of the
form $\tilde p_{R,1_Y}$ and $\tilde r_{1_X,S}$ are commutative while
Proposition 4.8 of \cite{ckww} says that these same squares satisfy 
the Beck condition. 
If $R$ and $S$ are also maps and $\alpha$ is invertible then
$\alpha\inv$ gives rise to another arrow of $\bG$, from $f$ to $u$
with reference to the square above, which may or may
not satisfy the Beck condition. The point here is that a
commutative square of maps gives rise to two, generally distinct,
Beck conditions.  It is well known that, for bicategories
of the form $\spn\E$ and $\rel\E$, all pullback
squares of maps satisfy both Beck conditions. A category
with finite products has automatically a number of pullbacks which
we might call {\em product-absolute} pullbacks because they are
preserved by all functors which preserve products. In \cite{ww} 
the Beck conditions for the product-absolute pullback squares of the form
$$\bfig
\Square(1000,0)/->`<-`<-`->/[A\x A`A\x A\x A`A`A\x A;d\x A`d`A\x d`d]
\efig$$
were investigated. 
(In fact, in this case it was shown that either Beck condition
implies the other.) The objects for which these conditions 
are met are called {\em Frobenius} objects.
\prp\label{mcifro} For a cartesian bicategory, the axiom {\em Maps are Comonadic} 
implies the axiom {\em Frobenius}.
\eth
\prf It suffices to show that the 2-cell $\delta_1$ below is invertible:
$$\bfig
\square(0,0)|alrb|/<-`->``<-/<750,500>%
[A`A\ox A`A\ox A`A\ox(A\ox A);d^*`d``1\ox d^*]
\morphism(750,500)|r|<0,-250>[A\ox A`(A\ox A)\ox A;d\ox 1]
\morphism(750,250)|r|<0,-250>[(A\ox A)\ox A`A\ox(A\ox A);a]
\morphism(175,250)|a|<150,0>[`;\delta_1]
\square(0,-500)|blrb|/<-`->`->`<-/<750,500>%
[A\ox A`A\ox(A\ox A)`A`A\ox A;1\ox d^*`r`r`d^*]
\morphism(250,-250)|a|<150,0>[`;\tilde r_{1_A,d^*}]
\efig$$
The paste composite of the squares is invertible (being 
essentially the identity 2-cell on $d^*$). The lower 2-cell is invertible 
by Proposition 4.7 of \cite{ckww} so that the whisker composite
$r\delta_1$ is invertible. Since $r$ is a map it reflects isomorphisms, by 
Maps are Comonadic, and hence $\delta_1$ is invertible.
\frp

\rmk\label{frobclofin} It was shown in \cite{ww} that, in a cartesian 
bicategory, the Frobenius objects are closed under finite products.
It follows that the full subbicategory of a cartesian bicategory
determined by the Frobenius objects is a cartesian bicategory
which satisfies the Frobenius axiom. 
\eth 

\section{Separable Objects and Discrete Objects 
in Cartesian Bicategories}\label{Sep}

In this section we look at separability for objects of cartesian bicategories.
Since for an object $A$ which is both separable and Frobenius, the hom-category
$\map\bB(X,A)$ is essentially discrete, for all $X$, we shall then be able to 
show that $\map\bB$ is essentially discrete by showing that all objects in $\bB$ 
are separable and Frobenius. But first we record the following direct argument:

\begin{proposition}\label{easy}
If $\bB$ is a bicategory in which all maps are comonadic and $\map\bB$ has
a terminal object, then $\map\bB$ is locally essentially discrete.
\end{proposition}

\prf
We must show that for all objects $X$ and $A$, the hom-category $\map\bB(X,A)$
is essentially discrete. As usual, we write $1$ for the terminal object of $\map\bB$
and $t_A:A\to 1$ for the essentially unique map, which by assumption is 
comonadic. Let $f,g:X\to A$ be maps from $X$ to $A$. If $\alpha:f\to g$ is any
2-cell, then $t_A\alpha$ is invertible, since $1$ is terminal in $\map\bB$. But
since $t_A$ is comonadic, it reflects isomorphisms, and so $\alpha$ is invertible.
Furthermore, if $\beta:f\to g$ is another 2-cell, then $t_A\alpha=t_A\beta$ by 
the universal property of $1$ once again, and now $\alpha=\beta$ since $t_A$
is faithful. Thus there is at most one 2-cell from $f$ to $g$, and any such 2-cell
is invertible. 
\frp

In any (bi)category with finite products the diagonal
arrows $d_A\f A\ra A\x A$ are (split) monomorphisms so that
in the bicategory $\bM$ the following square is a product-absolute 
pullback
$$\bfig
\square(0,0)[A`A`A`A\ox A;1_A`1_A`d_A`d_A]
\efig$$
that gives rise to a single $\bG$ arrow.

\dfn\label{sep}
An object $A$ in a cartesian bicategory is said to be
{\em separable} if the $\bG$ arrow above satisfies the Beck
condition. 
\eth
Of course the invertible mate condition here says precisely that the 
unit $\eta_{d_A}\f1_A\ra d_A^*d_A$ for the adjunction 
$d_A\laj d_A^*$ is invertible. Thus Axiom \ref{sepax}, as stated in
the Introduction, says that, for all $A$ in $\bB$, $\eta_{d_A}$ is
invertible. 

\rmk\label{sepcat} For a map $f$ it makes sense to define
{\em $f$ is fully faithful} to mean that $\eta_f$ is invertible. For
a {\em category $A$} the diagonal $d_A$ is fully faithful if and only
if $A$ is an ordered set.
\eth

\prp\label{sepmeans}
For an object $A$ in a cartesian bicategory, 
the following are equivalent:
\begin{enumerate}[$i)$]
\item $A$ is separable;
\item for all $f\f X\ra A$ in $\bM$, the diagram $f\la f\ra f$ is a product in
           $\bB(X,A)$;
\item $1_A\la 1_A\ra 1_A$ is a product in $\bB(A,A)$;
\item $1_A\ra\top_{A,A}$ is a monomorphism in $\bB(A,A)$;
\item for all $G\ra1_A$ in $\bB(A,A)$, the diagram $G\la G\ra1_A$ is 
a product in $\bB(A,A)$.
\end{enumerate}
\eth
\prf 

$[i)\Longrightarrow$ $ii)]$ A local product of maps is not generally a 
     map but here we have: 
$$f\wedge f\iso d_A^*(f\ox f)d_X\iso d_A^*(f\x f)d_X\iso%
 d_A^*d_A f\iso f$$ 

$[ii)\Longrightarrow$ $iii)]$ is trivial.

$[iii)\Longrightarrow$ $i)]$ Note the use of 
pseudo-functoriality of $\ox$:
$$d^*_Ad_A\iso d^*_A1_{A\ox A}d_A\iso d^*_A(1_A\ox1_A)d_A%
\iso1_A\wedge1_A\iso1_A$$

$[iii)\Longrightarrow$ $iv)]$ To say that $1_A\la 1_A\ra 1_A$ is a 
product in $\bB(A,A)$ is precisely to say that
$$\bfig
\square(0,0)[1_A`1_A`1_A`\top_{A,A};1_{1_A}`1_{1_A}``]
\efig$$
is a pullback in $\bB(A,A)$ which in turn is precisely to say that
$1_A\ra\top_{A,A}$ is a monomorphism in $\bB(A,A)$

$[iv)\Longrightarrow$ $v)]$ It is a generality that if an object
$S$ in a category is subterminal then for any $G\ra S$, necessarily 
unique, $G\la G\ra S$ is a product diagram.

$[v)\Longrightarrow$ $iii)]$ is trivial.
\frp
\cor\label{mcisep}{\rm [Of $iv)$]} For a cartesian bicategory, the axiom {\em Maps are
Comonadic} implies the axiom {\em Separable}.
\eth
\prf
We have $\top_{A,A}=t_A^*t_A$ for the map $t_A\f A\ra1$. It follows
that the unique $1_A\ra t_A^*t_A$ is $\eta_{t_A}$. Since $t_A$ is 
comonadic, $\eta_{t_A}$ is the equalizer shown:
$$1_A\to^{\eta_{t_A}}t_A^*t_A\two%
^{t_A^*t_A\eta_{t_A}}_{\eta_{t_A}t_A^*t_A}t_A^*t_At_A^*t_A$$
and hence a monomorphism.
\frp
\cor\label{copt}{\rm [Of $iv)$]} For separable $A$ in cartesian $\bB$, an
arrow $G\f A\ra A$ admits at most one copoint $G\ra 1_A$ depending upon
whether the unique arrow $G\ra\top_{A,A}$ factors through 
$1_A{\qqq\mon}\top_{A,A}$.
\frp
\eth

\prp\label{sepclofin} In a cartesian bicategory,
the separable objects are closed under finite products.
\eth
\prf
If $A$ and $B$ are separable objects then applying
the homomorphism
$\ox\f\bB\x\bB\ra\bB$
we have an adjunction $d_A\x d_B\laj d_A^*\ox d_B^*$ with
unit $\eta_{d_A}\ox\eta_{d_B}$ which being an isomorph of the
adjunction $d_{A\ox B}\laj d^*_{A\ox B}$ with unit
$\eta_{d_{A\ox B}}$ (via middle-four interchange) shows that the separable 
objects are closed under binary products. On the other hand, $d_I$ is an 
equivalence so that $I$ is also separable.
\frp

\cor\label{eos}
For a cartesian bicategory, the full subbicategory determined by
the separable objects is a cartesian bicategory which satisfies the axiom
{\em Separable}.
\frp
\eth

\prp\label{ordhom}
If $A$ is a separable object in a cartesian bicategory $\bB$,
then, for all $X$ in $\bB$, the hom-category $\bM(X,A)$ is an ordered
set, meaning that the category structure forms a reflexive, 
transitive relation.
\eth
\prf
Suppose that we have arrows $\alpha,\beta\f g{\s \two}f$ in
$\bM(X,A)$. In $\bB(X,A)$ we have
$$\bfig
\Atriangle(0,0)/->`->`/[g`f`f;\alpha`\beta`]
\morphism(500,500)|m|<0,-500>[g`f\wedge f;\gamma]
\morphism(500,0)|b|<-500,0>[f\wedge f`f;\pi]
\morphism(500,0)|b|<500,0>[f\wedge f`f;\rho]
\efig$$
By Proposition \ref{sepmeans} we can take $f\wedge f =f$ and
$\pi=1_f=\rho$ so that we have $\alpha=\gamma=\beta$. It follows 
that $\bM(X,A)$ is an ordered set. 
\frp
\dfn\label{discrete}
An object $A$ in a cartesian bicategory is said to be
{\em discrete} if it is both Frobenius and separable. We write
$\dis\bB$ for the full subbicategory of $\bB$ determined
by the discrete objects.
\eth
\begin{remark}
Beware that this is quite different to the notion of discreteness in a bicategory.
An object $A$ of a bicategory is discrete if each hom-category $\bB(X,A)$ is 
discrete; $A$ is essentially discrete if each $\bB(X,A)$ is equivalent to a discrete category. The notion of discreteness for cartesian bicategories defined above turns out to mean
that $A$ is essentially discrete in the bicategory $\map\bB$.
\end{remark}
From Proposition \ref{sepclofin} above and Proposition 3.4 of \cite{ww}
we immediately have
\prp\label{eod}
For a cartesian bicategory $\bB$, the full subbicategory
$\dis\bB$ of discrete objects is a cartesian bicategory in which
every object is discrete.
\frp
\eth
And from Proposition \ref{ordhom} above and Theorem 3.13 of \cite{ww}
we have
\prp\label{dishom}
If $A$ is a discrete object in a cartesian bicategory $\bB$
then, for all $X$ in $\bB$, the hom category $\bM(X,A)$ is an
equivalence relation.
\frp
\eth
If both the {\em Frobenius} axiom of \cite{ww}
and the {\em Separable} axiom of this paper hold for our cartesian
bicategory $\bB$, then every object of $\bB$ is discrete. In
this case, because $\bM$ is a bicategory, the equivalence relations
$\bM(X,A)$ are stable under composition from both sides. Thus
writing $|\bM(X,A)|$ for the set of objects of $\bM(X,A)$ we have
a mere category, $\E$ whose objects are those
of $\bM$ (and hence also those of $\bB$) and whose hom sets are
the quotients $|\bM(X,A)|/\bM(X,A)$. If the $\E(X,A)$
are regarded as discrete categories, so that $\E$ is a locally
discrete bicategory then the functors
$\bM(X,A)\ra |\bM(X,A)|/\bM(X,A)$ constitute the effect on homs functors
for an identity on objects biequivalence $\bM\ra \E$. To summarize
\thm\label{odibld} If a cartesian bicategory $\bB$ satisfies both the Frobenius
and Separable axioms then the bicategory of maps $\bM$ is
biequivalent to the locally discrete bicategory $\E$. \frp
\eth

In the following lemma we show that any copointed endomorphism of a discrete
object can be made into a comonad; later on, we shall see that this comonad 
structure is in fact unique.
\lem\label{diag} If $A$ is a discrete object in a
cartesian bicategory $\bB$ then, for any copointed endomorphism arrow
$\epsilon\f G\ra 1_A\f A\ra A$, there is a 2-cell 
$\delta=\delta_G\f G\ra GG$
satisfying
$$\bfig
\Atriangle(0,0)/->`->`/[G`G`G;1`1`]
\morphism(500,500)|m|<0,-500>[G`GG;\delta]
\morphism(500,0)|b|<-500,0>[GG`G;G\epsilon]
\morphism(500,0)|b|<500,0>[GG`G;\epsilon G]
\efig$$
and if both $G,H\f A\tra A$ are copointed, so that $GH\f A\ra A$
is also copointed, and $\phi\f G\ra H$ is any 2-cell, then the 
$\delta$'s satisfy
$$\bfig
\square(0,0)[G`H`GG`HH;\phi`\delta`\delta`\phi\phi]
\place(1000,0)[\mbox{and}]
\Atriangle(1500,0)/<-`->`->/%
[GHGH`GH`GH;\delta`(G\epsilon)(\epsilon H)`1]
\efig$$
\eth
\prf
We define $\delta=\delta_G$ to be the pasting composite
$$\bfig
\qtriangle(0,0)|amm|[A`AA`AAA;d`d_3`1d]
\square(500,0)|amma|[AA`AA`AAA`AAA;G1`1d`1d`G11]
\square(1000,0)|amma|[AA`A`AAA`AA;d^*`1d`d`d^*1]
\qtriangle(500,-500)|abr|[AAA`AAA`AAA;G11`GGG`11G]
\square(1000,-500)|arma|[AAA`AA`AAA`AA;d^*1`11G`1G`d^*1]
\qtriangle(1000,-1000)|amm|[AAA`AA`A;d^*1`d_3^*`d^*]
\morphism(0,500)|b|/{@{->}@/_4em/}/<1500,-1500>[A`A;G\wedge G\wedge G]
\morphism(0,500)|b|/{@{->}@/_8em/}/<1500,-1500>[A`A;G]
\morphism(0,500)|a|/{@{->}@/^3em/}/<1500,0>[A`A;G]
\morphism(1500,500)|r|/{@{->}@/^3em/}/<0,-1500>[A`A;G]
\morphism(800,-250)|m|<150,150>[`;G\epsilon G]
\morphism(300,-800)|m|<150,150>[`;\delta_3]
\place(1250,250)[1]
\place(750,675)[2]
\place(1700,-250)[3]
\place(750,250)[4]
\place(1250,-250)[5]
\place(600,-450)[6]
\efig$$
wherein $\ox$ has been abbreviated by juxtaposition and all
subregions not explicitly inhabited by a 2-cell are deemed to
be inhabited by the obvious invertible 2-cell. A reference number has been
assigned to those invertible 2-cells which arise from the hypotheses.
As in \cite{ww}, $d_3$'s denote 3-fold diagonal maps and, similarly,
we write $\delta_3$ for a local 3-fold diagonal. 

The invertible 2-cell labelled by `1' is that defining $A$ to be
Frobenius. The 3-fold composite of arrows in the region labelled
by `2' is $G\wedge1_A$ and, similarly, in that labelled by `3' we have
$1_A\wedge G$. Each of these is isomorphic to $G$ because $A$ is
separable and $G$ is copointed. The isomorphisms in `4' and `5'
express the pseudo-functoriality of $\ox$ in the cartesian bicategory
$\bB$. Finally `6' expresses the ternary local product in terms of
the ternary $\ox$ as in \cite{ww}. Demonstration of the equations is
effected easily by pasting composition calculations.
\frp
\thm\label{wedge=.} If $G$ and $H$ are copointed endomorphisms
on a discrete $A$ in a cartesian $\bB$ then
$$G\toleft^{G\epsilon}GH\to^{\epsilon H}H$$
is a product diagram in $\bB(A,A)$.
\eth
\prf
If we are given $\alpha\f K\ra G$ and $\beta\f K\ra H$ then $K$ is
also copointed and we have
$$K\to^\delta KK\to^{\alpha\beta}GH$$
as a candidate pairing. That this candidate satisfies the universal
property follows from the equations of Lemma \ref{diag} which
are precisely those in the equational description of binary products.
We remark that the `naturality' equations for the projections follow
immediately from uniqueness of copoints.
\frp
\cor\label{comsim} If $A$ is discrete in a cartesian $\bB$, then an endo-arrow
$G\f A\ra A$ admits a comonad structure if and only if $G$ has
the copointed property, and any such comonad structure is unique.
\eth
\prf
The Theorem shows that the arrow $\delta\f G\ra GG$ constructed in Lemma
\ref{diag} is the product diagonal on $G$ in the category $\bB(A,A)$ and,
given $\epsilon\f G\ra 1_A$, this is the only comonad comultiplication on $G$.   
\frp 
\rmk
It is clear that $1_A$ is terminal with respect to the copointed
objects in $\bB(A,A)$.
\eth
\prp\label{subterm} If an object $B$ in a bicategory $\bB$ has $1_B$ 
subterminal in $\bB(B,B)$ then, for any map $f\f A\ra B$, $f$ is subterminal in $\bB(A,B)$ and $f^*$ is subterminal in $\bB(B,A)$. In particular, in a 
cartesian bicategory in which every object is separable, every adjoint
arrow is subterminal.
\eth
\prf
Precomposition with a map preserves terminal objects and monomorphisms,
as does postcomposition with a right adjoint.
\frp  
\section{Bicategories of Comonads}\label{Coms}
The starting point of this section is the observation, made in the introduction, 
that a comonad in the bicategory $\spn\E$ is 
precisely a span of the form 
$$A \xleftarrow{g} X \xrightarrow{g} A$$
in which both legs are equal.
 
We will write $\bC=\com\bB$ for the bicategory of comonads in $\bB$,
$\com$ being one of the duals of Street's construction $\mathrm{Mnd}$
in \cite{ftm}. Thus $\bC$ has objects given by the comonads $(A,G)$ 
of $\bB$. The structure 2-cells for comonads will be denoted
$\epsilon=\epsilon_G$, for counit and $\delta=\delta_G$, for 
comultiplication. An arrow 
in $\bC$ from $(A,G)$ to $(B,H)$ is a pair $(F,\phi)$ as shown in
$$\bfig
\square(0,0)[A`B`A`B;F`G`H`F]
\morphism(125,250)|a|<250,0>[`;\phi]
\efig$$
satisfying
\begin{equation}\label{comarrow}
\bfig
\square(0,0)[FG`HF`F1_A`1_BF;\phi`F\epsilon`\epsilon F`=]
\place(1000,250)[\mbox{and}]
\square(1500,0)/`->``->/[FG``FGG`HFG;`F\delta``\phi G]
\square(2000,0)/``->`->/[`HF`HFG`HHF;``\delta F`H\phi]
\morphism(1500,500)<1000,0>[FG`HF;\phi]
\efig
\end{equation}
(where, as often, we have suppressed the associativity constraints of 
our normal, cartesian, bicategory $\bB$).
A 2-cell $\tau\f(F,\phi)\ra(F',\phi')\f(A,G)\ra(B,H)$ in $\bC$
is a 2-cell $\tau\f F\ra F'$ in $\bB$ satisfying
\begin{equation}\label{comtrans}
\bfig
\square(0,0)[FG`HF`F'G`HF';\phi`\tau G`H\tau`\phi']
\efig
\end{equation}

There is a pseudofunctor $I\f\bB\ra\bC$ given by
$$I(\tau\f F\ra F'\f A\ra B)=%
\tau\f (F,1_{F})\ra (F',1_{F'})\f (A,1_A)\ra (B,1_B)$$
From \cite{ftm} it is well known that a bicategory $\bB$
has Eilenberg-Moore objects for comonads
if and only if $I\f\bB\ra\bC$ has a right biadjoint, which we will
denote by $E\f\bC\ra\bB$. We write $E(A,G)=A_G$ and the counit
for $I\laj E$ is denoted by 
$$\bfig
\square(0,0)[A_G`A`A_G`A;g_G`1_{A_G}`G`g_G]
\morphism(125,250)|a|<250,0>[`;\gamma_G]
\Ctriangle(2500,0)/<-`->`->/<500,250>[A`A_G`A;g_G`G`g_G]
\place(1500,250)[\mbox{ or, using normality of $\bB$, better by }]
\morphism(2700,250)|m|<200,0>[`;\gamma_G]
\efig$$
with $(g_G,\gamma_G)$ abbreviated to $(g,\gamma)$ when there is
no danger of confusion. It is standard that each $g=g_G$ is necessarily
a map (whence our lower case notation) and the mate $gg^*\ra G$ of 
$\gamma$ is an isomorphism which identifies $\epsilon_g$ and
$\epsilon_G$.

We will write $\bD$ for the locally full subbicategory
of $\bC$ determined by all the objects and those arrows of the form $(f,\phi)$,
where $f$ is a map, and write $j\f\bD\ra\bC$ for the inclusion. 
It is clear that the
pseudofunctor $I\f\bB\ra\bC$ restricts to give a pseudofunctor
$J\f\bM\ra\bD$. We say that the bicategory $\bB$ {\em has Eilenberg-Moore 
objects for comonads, as seen by $\bM$}, if $J\f\bM\ra\bD$ has a right
biadjoint. (In general, this property does not follow from that of
Eilenberg-Moore objects for comonads.) 

\begin{remark}\label{rmk:D-for-dummies}
In the case $\bB=\spn\E$, a comonad in $\bB$ can, as we have seen, be 
identified with a morphism in \E. This can be made into the object part of a 
biequivalence between the bicategory $\bD$ and the category $\E^\arr$
of arrows in $\E$. If we further identify $\bM$ with $\E$, then the 
inclusion $j:\bD\to\bC$ becomes the diagonal $\E\to\E^\arr$; of course this
does have a right adjoint, given by the domain functor.
\end{remark}

\thm\label{simcom} If $\bB$ is a cartesian bicategory in which every object
is discrete, the bicategory $\bD=\bD(\bB)$ admits the following simpler
description:
\begin{enumerate}
\item[$i)$] An object is a pair $(A,G)$ where $A$ is an object of $\bB$
and $G\f A\ra A$ admits a copoint;
\item[$ii)$] An arrow $(f,\phi)\f(A,G)\ra(B,H)$ is a map $f\f A\ra B$ and
a 2-cell $\phi\f fG\ra Hf$;
\item[$iii)$] A 2-cell $\tau\f(f,\phi)\ra(f',\phi')\f(A,G)\ra(B,H)$ is
a 2-cell satisfying $\tau\f f\ra f'$ satisfying equation (\ref{comtrans}).
\end{enumerate} 
\eth
\prf
We have i) by Corollary \ref{comsim} while iii) is precisely the description
of a 2-cell in $\bD$, modulo the description of the domain and codomain
arrows. So, we have only to show ii), which is to show that the equations
(\ref{comarrow}) hold automatically under the hypotheses. For the first
equation of (\ref{comarrow}) we have uniqueness of any 2-cell $fG\ra f$ 
because $f$ is subterminal by Proposition \ref{subterm}. For the second,
observe that the terminating vertex, $HHf$, is the product $Hf\wedge Hf$
in $\bM(A,B)$ because $HH$ is the product $H\wedge H$ in $\bM(B,B)$ by
Theorem~\ref{wedge=.} and precomposition with a map preserves all limits.
For $HHf$ seen as a product, the projections are, again 
by Theorem~\ref{wedge=.}, $H\epsilon f$ and $\epsilon Hf$. Thus, it 
suffices to show that the diagram for the second equation commutes when
composed with both $H\epsilon f$ and $\epsilon Hf$. We have
$$\bfig
\square(0,0)/`->``->/[fG``fGG`HfG;`f\delta``\phi G]
\square(500,0)/``->`->/[`Hf`HfG`HHf;``\delta f`H\phi]
\morphism(0,500)<1000,0>[fG`Hf;\phi]
\qtriangle(500,-500)|blr|[HfG`HHf`Hf;H\phi`Hf\epsilon`H\epsilon f]
\morphism(0,-500)|b|<1000,0>[fG`Hf;\phi]
\morphism(0,0)<0,-500>[fGG`fG;fG\epsilon]
\square(2000,0)/`->``->/[fG``fGG`HfG;`f\delta``\phi G]
\square(2500,0)/``->`->/[`Hf`HfG`HHf;``\delta f`H\phi]
\morphism(2000,500)<1000,0>[fG`Hf;\phi]
\ptriangle(2000,-500)|blr|[fGG`HfG`fG;\phi G`f\epsilon G`\epsilon fG]
\morphism(2000,-500)|b|<1000,0>[fG`Hf;\phi]
\morphism(3000,0)|r|<0,-500>[HHf`Hf;\epsilon Hf]
\efig$$
in which each of the lower triangles commutes by the first equation of 
(\ref{comarrow}) already established. Using comonad equations for
$G$ and $H$, it is obvious that each composite is $\phi$.
\frp
Finally, let us note that $\bD$ is 
a subbicategory, neither full nor locally full, of the Grothendieck 
bicategory $\bG$ and write $K\f\bD\ra\bG$ for the inclusion. We also 
write $\iota\f\bM\ra\bG$ 
for the composite pseudofunctor $KJ$. Summarizing, we have introduced
the following commutative diagram of bicategories and pseudofunctors 
$$\bfig
\square(0,0)[\bM`\bD`\bB`\bC;J`i`j`I]
\morphism(500,500)<500,0>[\bD`\bG;K]
\morphism(0,500)|a|/{@{->}@/^2em/}/<1000,0>[\bM`\bG;\iota]
\efig ;$$
note also that in our main case of interest $\bB=\spn\E$, each of $\bM$,
$\bD$, and $\bG$ is biequivalent to a mere category.
Ultimately, we are interested in having a right biadjoint,
say $\tau$, of $\iota$. For such a biadjunction $\iota\laj\tau$
the counit at an object $R\f X\ra A$ in $\bG$ will take the form
\begin{equation}\label{tabcounit}
\bfig
\Ctriangle/<-`->`->/<500,250>[X`\tau R`A;u_R`R`v_R]
\morphism(200,250)|m|<200,0>[`;\omega_R]
\efig
\end{equation}
(where, as for a biadjunction $I\laj E\f\bC\ra\bB$, a
triangle rather than a square can be taken as the boundary of
the 2-cell by the normality of $\bB$). In fact, we are interested
in the case where we have $\iota\laj\tau$ and moreover the counit
components $\omega_R\f v_R\ra Ru_R$ enjoy the property that their 
mates $v_Ru^*_R\ra R$ with respect to the adjunction $u_R\laj u^*_R$
are invertible. In this way we represent a general arrow of $\bB$
in terms of a span of maps. Since biadjunctions compose we will 
consider adjunctions $J\laj F$ and $K\laj G$ and we begin with the 
second of these.

\thm\label{G(R)}
For a cartesian bicategory $\bB$ in which every object is discrete,
there is an adjunction $K\laj G\f\bG\ra\bD$ where, for
$R\f X\ra A$ in $\bG$, the comonad $G(R)$ and its witnessing copoint
$\epsilon\f G(R)\ra 1_{XA}$ are given by the left diagram below and the
counit $\mu\f KG(R)\ra R$ is given by the right diagram below,
all in notation suppressing $\ox$:
$$\bfig
\ptriangle(-1500,1000)/->`->`<-/[XA`XA`XXA;1_{XA}`dA`p_{1,3}]
\morphism(-1400,1350)|a|<150,0>[`;\simeq]
\morphism(-1500,1000)|l|<0,-500>[XXA`XAA;XRA]
\morphism(-1000,1500)|r|<0,-1500>[XA`XA;1_{XA}]
\morphism(-1375,750)|m|<250,0>[`;\tilde p_{1,3}]
\btriangle(-1500,0)[XAA`XA`XA;Xd^*`p_{1,3}`1_{XA}]
\morphism(-1400,150)|a|<150,0>[`;]
\morphism(-1500,1500)|l|/{@{->}@/_3.5em/}/<0,-1500>[XA`XA;G(R)]
\ptriangle(0,1000)/->`->`<-/[XA`X`XXA;p`dA`p_2]
\morphism(0,1000)|l|<0,-500>[XXA`XAA;XRA]
\morphism(500,1500)|r|<0,-1500>[X`A;R]
\btriangle(0,0)[XAA`XA`A;Xd^*`p_2`r]
\morphism(100,1350)|a|<150,0>[`;\simeq]
\morphism(125,750)|m|<250,0>[`;\tilde p_2]
\morphism(100,150)|b|<150,0>[`;]
\morphism(0,1500)|l|/{@{->}@/_3.5em/}/<0,-1500>[XA`XA;G(R)]
\efig$$
Moreover, the mate $rG(R)p^*\ra R$ of the counit $\mu$ is invertible.

In the left diagram, the $p_{1,3}$ collectively denote projection from the 
three-fold product in $\bG$ to the product of the first and third factors.  
In the right diagram, the $p_2$ collectively denote projection from the 
three-fold product in $\bG$ to the second factor. The upper triangles of the 
two diagrams are the canonical isomorphisms. The lower left triangle 
is the mate of the canonical isomorphism $1{\s \to^{\simeq}}p_{1,3}(Xd)$.
The lower right triangle is the mate of the canonical isomorphism
$r{\s \to^{\simeq}} p_2(Xd)$.
\eth
\prf Given a comonad $H\f T\ra T$ and an arrow
$$\bfig
\square(0,0)[T`X`T`A;x`H`R`a]
\morphism(125,250)|a|<250,0>[`;\psi]
\efig$$
in $\bG$, we verify the adjunction claim by showing that there is 
a unique arrow
$$\bfig
\square(0,0)[T`XA`T`XA;f`H`G(R)`f]
\morphism(125,250)|a|<250,0>[`;\phi]
\efig$$
in $\bD$, whose composite with the putative counit $\mu$ is
$(x,\psi,a)$. It is immediately clear that the unique solution
for $f$ is $(x,a)$ and to give $\phi\f(x,a)H\ra Xd^*(XRA)dA(x,a)$
is to give the mate $Xd(x,a)H\ra (XRA)dA(x,a)$ which is
$(x,a,a)H\ra (XRA)(x,x,a)$ and can be seen as
a $\bG$ arrow:
$$\bfig
\square(0,0)[T`XXA`T`XAA;(x,x,a)`H`XRA`(x,a,a)]
\morphism(125,250)|a|<250,0>[`;(\alpha,\beta,\gamma)]
\efig$$
where we exploit the description of products in $\bG$. From this 
description it is clear, since $\tilde p_2(\alpha,\beta,\gamma)=\beta$ 
as a composite in $\bG$, that the unique solution for $\beta$ is 
$\psi$. We have seen in Theorem \ref{comsim} that the conditions 
(\ref{comarrow}) hold automatically in $\bD$ under the assumptions 
of the Theorem. From the first of these we have:
$$\bfig
\square(0,0)|almb|[T`XXA`T`XAA;(x,x,a)`H`XRA`(x,a,a)]
\morphism(125,250)|a|<250,0>[`;(\alpha,\beta,\gamma)]
\square(500,0)|amrb|[XXA`XA`XAA`XA;p_{1,3}`XRA`1_{XA}`p_{1,3}]
\morphism(625,250)|a|<250,0>[`;p_{1,3}]
\place(1375,250)[=]
\square(2000,0)|arrb|[T`XA`T`XA;(x,a)`1_T`1_{XA}`(x,a)]
\morphism(2125,250)|a|<250,0>[`;\kappa_{(x,a)}]
\morphism(2000,500)|l|/{@{->}@/_3em/}/<0,-500>[T`T;H]
\morphism(1750,250)|a|<200,0>[`;\epsilon_H]
\efig$$
So, with a mild abuse of notation, we have
$(\alpha,\gamma)=(1_x\epsilon_H, 1_a\epsilon_H)$, uniquely,
and thus the unique solutions for $\alpha$ and $\gamma$ are
$1_x\epsilon_H$ and $1_a\epsilon_H$ respectively. This shows that
$\phi$ is necessarily the mate under the adjunctions considered of
$(1_x\epsilon_H,\psi,1_a\epsilon_H,)$. Since $\bD$ and $\bG$
are essentially locally discrete this suffices to
complete the claim that $K\laj G$. 
It only remains to show that the mate $rG(R)p^*\ra R$ of the counit 
$\mu$ is invertible. In the three middle squares of the diagram
$$\bfig
\morphism(0,500)|l|/{@{->}@/_3.5em/}/<0,-1500>[XA`XA;G(R)]
\square(0,0)|alrm|/<-`->`->`<-/[XA`X`XXA`XX;p^*`dA`d`p^*]
\morphism(125,250)|a|<250,0>[`;\tilde p^*_{d,1_A}]
\square(0,-500)|mlmm|/<-`->`->`<-/[XXA`XX`XAA`XA;p^*`XRA`XR`p^*]
\morphism(125,-250)|a|<250,0>[`;\tilde p^*_{XR,1_A}]
\square(0,-1000)|mlrb|/<-`->`->`->/[XAA`XA`XA`A;p^*`Xd^*`r`r]
\morphism(175,-750)|a|<150,0>[`;\simeq]
\square(500,-500)|amrb|[XX`X`XA`A;r`XR`R`r]
\morphism(625,-250)|a|<250,0>[`;r_{1_X,R}]
\morphism(500,500)|r|<500,-500>[X`X;1_X]
\morphism(1000,-500)|r|<-500,-500>[A`A;1_A]
\efig$$
the top two are invertible 2-cells by
Proposition 4.18 of \cite{ckww} while the lower one is the obvious invertible
2-cell constructed from $Xd^*p^*\iso1_{X,A}$. The right square is an
invertible 2-cell by Proposition 4.17 of \cite{ckww}. This shows that the
mate $rG(R)p^*\ra R$ of $\mu$ is
invertible.
\frp

\rmk \label{unit} It now follows that the unit of the adjunction $K\laj G$
is given (in notation suppressing $\ox$) by:
$$\bfig
\qtriangle(1500,1000)[T`TT`TTT;d`d_3`dT]
\dtriangle(1500,0)/<-`->`->/[TTT`T`TT;d_3`Td^*`d]
\morphism(1500,1500)|l|<0,-1500>[T`T;H]
\morphism(2000,1000)|l|<0,-500>[TTT`TTT;HHH]
\morphism(2000,1000)|r|/{@{->}@/^3em/}/<0,-500>[TTT`TTT;THT]
\morphism(1750,1350)|a|<150,0>[`;\simeq]
\morphism(1550,750)|m|<175,0>[`;\tilde d_3]
\morphism(2050,750)|m|<225,0>[`;\epsilon H\epsilon]
\morphism(1750,150)|a|<150,0>[`;\simeq]
\efig$$
where the $d_3$ collectively denote 3-fold diagonalization $(1,1,1)$ in $\bG$.
The top triangle is a canonical isomorphism while the lower triangle
is the mate of the canonical isomorphism $(T\ox d)d{\s \to^{\simeq}}d_3$
and is itself invertible, by separability of $T$.
\eth
Before turning to the question of an adjunction $J\laj F$, we note:
\lem\label{maplikeanm}
In a cartesian bicategory in which Maps are Comonadic,
if $gF\iso h$ with $g$ and $h$ maps, then $F$ is also a map.
\eth
\prf
By Theorem 3.11 of \cite{ww} it suffices to show that $F$ is
a comonoid homomorphism, which is to show that the
canonical 2-cells $\tilde t_F\f tF\ra t$ and $\tilde d_F\f dF\ra(F\ox F)d$
are invertible. For the first we have:
$$tF\iso tgF\iso th\iso t$$
Simple diagrams show that we do get the right isomorphism
in this case and also for the next:
$$(g\ox g)(dF)\iso dgF\iso dh\iso(h\ox h)d\iso (g\ox g)(F\ox F)d$$
which gives $dF\iso(F\ox F)d$ since the map $g\ox g$ reflects
isomorphisms.
\frp
\thm\label{emasbm} If $\bB$ is a cartesian bicategory which has 
Eilenberg-Moore objects for Comonads and for which Maps are Comonadic 
then $\bB$ has Eilenberg-Moore objects for Comonads as Seen by $\bM$, 
which is to say that $J\f\bM\ra\bD$ has a right adjoint. Moreover, the counit
for the adjunction, say $JF\ra 1_{\bD}$, necessarily having components
of the form $\gamma\f g\ra Gg$ with $g$ a map, has $gg^*\ra G$ invertible.
\eth
\prf
It suffices to show that the adjunction $I\laj E\f\bC\ra\bB$
restricts to $J\laj F\f\bD\ra\bM$. For this it suffices to show that,
given $(h,\theta)\f JT\ra(A,G)$, the $F\f T\ra A_G$ with $gF\iso h$
which can be found using $I\laj E$ has $F$ a map. This follows from
Lemma \ref{maplikeanm}.
\frp
\thm\label{tabulation} A cartesian bicategory which has Eilenberg-Moore 
objects for Comonads and for which Maps are Comonadic has tabulation in 
the sense that the inclusion $\iota\f\bM\ra\bG$ has a right adjoint $\tau$ 
and the counit components $\omega_R\f v_R\ra Ru_R$ as in (\ref{tabcounit}) 
have the property that the mates $v_Ru^*_R\ra R$, with respect to the 
adjunctions $u_R\laj u^*_R$, are invertible. 
\eth
\prf
Using Theorems \ref{G(R)} and \ref{emasbm} we can construct the adjunction 
$\iota\laj\tau$ by composing $J\laj F$ with $K\laj G$. Moreover, the
counit for $\iota\laj\tau$ is the pasting composite:
$$\bfig 
\Ctriangle(0,0)|lmb|/<-`->`->/<500,250>[X\ox A`T`X\ox A;(u,v)`G(R)`(u,v)]
\morphism(0,250)|l|/{@{->}@/^4.0em/}/<1000,250>[T`X;u]
\morphism(0,250)|l|/{@{->}@/_4.0em/}/<1000,-250>[T`A;v]
\morphism(300,160)|m|<0,180>[`;\gamma]
\square(500,0)|amrb|<500,500>[X\ox A`X`X\ox A`A;p`G(R)`R`r]
\morphism(650,250)|m|<200,0>[`;\mu]
\efig$$
where the square is the counit for $K\laj G$; and the triangle, the counit
for $J\laj F$, is an Eilenberg-Moore coalgebra for the comonad $G(R)$. 
The arrow component of the Eilenberg-Moore coalgebra is necessarily of
the form $(u,v)$, where $u$ and $v$ are maps, and it also follows that 
we have $(u,v)(u,v)^*\iso G(R)$. Thus we have
$$vu^*\iso  r(u,v)(p(u,v))^*\iso r(u,v)(u,v)^*p^*\iso rG(R)p^*\iso R$$
where the first two isomorphisms are trivial, the third arises from the 
invertibility of the mate of $\gamma$ as an Eilenberg-Moore structure,
and the fourth is invertibility of $\mu$, as in Theorem \ref{G(R)}. 
\frp
\thm\label{mapbhaspb}
For a cartesian bicategory $\bB$ with Eilenberg-Moore objects for 
Comonads and for which Maps are Comonadic, $\map\bB$ has pullbacks 
satisfying the Beck condition (meaning that for a pullback square
\begin{equation}\label{beckforpb}
\bfig
\square(0,0)[P`M`N`A;r`p`b`a]
\morphism(200,250)<100,0>[`;\simeq]
\efig
\end{equation}
the mate $pr^*\ra a^*b$ of $ap\iso br$ in $\bB$, 
with respect to the adjunctions $r\laj r^*$ and $a\laj a^*$, 
is invertible).
\eth
\prf
Given the cospan $a\f N\ra A\la M\f b$ in $\map\bB$, let
$P$ together with $(r,\sigma,p)$ be a tabulation 
for $a^*b\f M\ra N$. Then $pr^*\ra a^*b$, the mate of
$\sigma\f p\ra a^*br$ with respect to $r\laj r^*$, is invertible 
by Theorem \ref{tabulation}. We have also $ap\ra br$, the mate of 
$\sigma\f p\ra a^*br$ with respect to $a\laj a^*$. Since $A$ 
is discrete, $ap\ra br$ is also invertible and is the only 2-cell 
between the composite maps in question. If we have also 
$u\f N\la T\ra M\f v$, for maps $u$ and $v$ with $au\iso bv$, then the 
mate $u\ra a^*bv$ ensures that the span $u\f N\la T\ra M\f v$ factors 
through $P$ by an essentially unique map $w\f T\ra P$ with
$pw\iso u$ and $rw\iso v$.
\frp
\prp\label{tabonadic} In a cartesian bicategory with 
Eilenberg-Moore objects for Comonads and for which Maps are Comonadic,
every span of maps $x\f X\la S\ra A\f a$ gives rise to the following
tabulation diagram:
$$\bfig
\Ctriangle(0,0)|lmb|/<-`->`->/<500,250>[X`S`A;x`ax^*`a]
\morphism(300,160)|m|<0,180>[`;a\eta_x]
\efig$$
\eth
\prf A general tabulation counit $\omega_R\f v_R\ra Ru_R$ is given in terms
of the Eilenberg-Moore coalgebra for the comonad $(u,v)(u,v)^*$ and 
necessarily $(u,v)(u,v)^*\iso G(R)$. It follows that for $R=ax^*$, it 
suffices to show that $G(ax^*)\iso (x,a)(x,a)^*$.
Consider the diagram (with $\ox$ suppressed):
$$\bfig
\Atriangle(0,0)|bba|/->`->`/[XSA`XXA`XAA;XxA`XaA`]
\Vtriangle(0,500)|mmm|/`->`->/[SA`XS`XSA;`(x,S)A`X(S,a)]
\Atriangle(0,1000)|lrm|/->`->`/[S`SA`XS;(S,a)`(x,S)`]
\Ctriangle(-500,0)|lml|/->``->/[SA`XA`XXA;xA``dA]
\Dtriangle(1000,0)|mrr|/`->`->/[XS`XA`XAA;`Xa`Xd]
\morphism(500,1500)|l|/{@{->}@/_3em/}/<-1000,-1000>[S`XA;(x,a)]
\morphism(500,1500)|r|/{@{->}@/^3em/}/<1000,-1000>[S`XA;(x,a)]
\efig$$
The comonoid $G(ax^*)$ can be read, from left to right, along the
`W' shape of the lower edge as $G(ax^*)\iso Xd^*.XaA.Xx^*A.dA$.
But each of the squares in the diagram is a (product-absolute) pullback 
so that with Proposition \ref{mapbhaspb} at hand we can continue:
$$Xd^*.XaA.Xx^*A.dA\iso Xa.X(S,a)^*.(x,S)A.x^*A\iso%
Xa.(x,S).(S,a)^*. x^*A\iso (x,a)(x,a)^*$$
as required.
\frp

\section{Characterization of Bicategories of Spans}\label{charspan}
\subsection{}\label{Cfun}
If $\bB$ is a cartesian bicategory with $\map\bB$ essentially locally
discrete then each slice $\map\bB/(X\ox A)$ 
is also essentially locally discrete and we can write 
$\spn\map\bB(X,A)$ 
for the categories obtained by taking the quotients of the equivalence 
relations comprising the hom categories of the $\map\bB/(X\ox A)$. 
Then we can construct
functors $C_{X,A}\f\spn\map\bB(X,A)\ra\bB(X,A)$, where for
an arrow in $\spn\map\bB(X,A)$ as shown,
$$\bfig
\Ctriangle(0,0)|lml|/->`->`<-/[M`A`N;a`h`b]
\Dtriangle(500,0)|mrr|/->`->`<-/[M`X`N;h`x`y]
\efig$$
we define $C(y,N,b)=by^*$ and 
$C(h)\f ax^*=(bh)(yh)^*\iso bhh^*y^*\to^{b\epsilon_hy^*} by^*$.
If $\map\bB$ is known to have pullbacks then the 
$\spn\map\bB(X,A)$ become the hom-categories for a
bicategory $\spn\map\bB$ and we can consider whether the
$C_{X,A}$ provide the effects on homs for an identity-on-objects
pseudofunctor $C\f\spn\map\bB\ra\bB$. Consider
\begin{equation}\label{beck}
\bfig
\Atriangle(0,0)/->`->`/[N`Y`A;y`b`]
\Vtriangle(500,0)/`->`->/[N`M`A;``]
\Atriangle(1000,0)/->`->`/[M`A`X;a`x`]
\Atriangle(500,500)/->`->`/[P`N`M;p`r`]
\efig
\end{equation}
where the square is a pullback. In somewhat
abbreviated notation, what is needed further are coherent,
invertible 2-cells $\widetilde C\f CN.CM\ra C(NM)=CP$, 
for each composable pair of spans $M$, $N$, and
coherent, invertible 2-cells $C^\circ\f 1_A\ra C(1_A)$,
for each object $A$.
Since the identity span on $A$ is $(1_A,A,1_A)$,
and $C(1_A)=1_A.1^*_A\iso1_A.1_A\iso1_A$ we take
the inverse of this composite for $C^\circ$. To give the
$\widetilde C$ though is to give 2-cells
$yb^*ax^*\ra ypr^*x^*$ and since spans of the form
$(1_N,N,b)$ and $(a,M,1_M)$ arise as special cases, it is
easy to verify that to give the $\widetilde C$ it is
necessary and sufficient to give coherent, invertible 2-cells
$b^*a\ra pr^*$ for each pullback square in $\map\bB$. The inverse
of such a 2-cell $pr^*\ra b^*a$ is the mate of a 2-cell
$bp\ra aq$. But by discreteness a 2-cell
$bp\ra aq$ must be essentially an identity. Thus, definability
of $\widetilde C$ is equivalent to the invertibility in $\bB$
of the mate $pr^*\ra b^*a$ of the identity $bp\ra ar$, for each 
pullback square as displayed in (\ref{beck}). In short,
if $\map\bB$ has pullbacks and these satisfy the Beck condition
as in Proposition \ref{mapbhaspb} 
then we have a canonical pseudofunctor $C\f\spn\map\bB\ra\bB$.

\thm\label{spanmain}
For a bicategory $\bB$ the following are equivalent:
\begin{enumerate}[$i)$]
\item There is a biequivalence $\bB\simeq\spn\E$, 
           for $\E$ a category with finite limits;
\item The bicategory $\bB$ is cartesian, each comonad has an
           Eilenberg-Moore object, and every map is comonadic.
\item The bicategory $\map\bB$ is an essentially locally
           discrete bicategory with finite limits, satisfying in $\bB$
	   the Beck condition for pullbacks of maps, and the 
	   canonical $$C\f\spn\map\bB\ra\bB$$ is a 
	   biequivalence of bicategories.
\end{enumerate}
\eth
\prf
That $i)$ implies $ii)$ follows from our discussion in
the Introduction.  That $iii)$ implies $i)$
is trivial so we show that $ii)$ implies $iii)$.

We have already observed in Theorem \ref{odibld} that,
for $\bB$ cartesian with every object discrete, 
$\map\bB$ is essentially locally discrete and we have seen
by Propositions \ref{mcifro} and \ref{mcisep} that, 
in a cartesian bicategory in which Maps are Comonadic, every
object is discrete.
In Theorem~\ref{mapbhaspb} we have seen that, for $\bB$ satisfying
the conditions of $ii)$, $\map\bB$ has pullbacks, and hence all
finite limits and, in $\bB$ the Beck condition holds for pullbacks. 
Therefore we have the canonical pseudofunctor $C\f\spn\map\bB\ra\bB$ 
developed in \ref{Cfun}. To complete the proof it suffices to show that
the $C_{X,A}\f\spn\map\bB(X,A)\ra\bB(X,A)$ are equivalences of categories.

Define functors $F_{X,A}\f\bB(X,A)\ra\spn\map\bB(X,A)$ 
by $F(R)=F_{X,A}(R)=(u,\tau R,v)$ where
$$\bfig
\Ctriangle/<-`->`->/<500,250>[X`\tau R`A;u`R`v]
\morphism(200,250)|m|<200,0>[`;\omega]
\efig$$
is the $R$-component of the counit for 
$\iota\laj\tau\f\bG\ra\map\bB$. For a 2-cell
$\alpha\f R\ra R'$ we define $F(\alpha)$ to be
the essentially unique map satisfying
$$\bfig
\morphism(-500,500)|m|[\tau R`\tau R';F(\alpha)]
\Ctriangle(0,0)|rrr|/<-`->`->/[X`\tau R'`A;u'`R'`v']
\morphism(175,400)|a|<150,0>[`;\omega']
\morphism(-500,500)|l|<1000,500>[\tau R`X;u]
\morphism(-500,500)|l|<1000,-500>[\tau R`A;v]
\place(750,500)[=]
\Ctriangle(1000,0)/<-``->/[X`\tau R`A;u``v]
\morphism(1500,1000)|l|/{@{->}@/_1em/}/<0,-1000>[X`A;R]
\morphism(1500,1000)|r|/{@{->}@/^1.5em/}/<0,-1000>[X`A;R']
\morphism(1150,400)|a|<150,0>[`;\omega]
\morphism(1450,400)|a|<150,0>[`;\alpha]
\efig$$ 
(We remark that essential uniqueness here means that
$F(\alpha)$ is determined to within unique invertible 2-cell.)
Since $\omega\f v\ra Ru$ has mate $vu^*\ra R$ invertible, 
because $(v,\tau R,u)$ is a tabulation of $R$,
it follows that we have a natural isomorphism
$CFR\ra R$. On the other hand, starting with a span
$(x,S,a)$ from $X$ to $A$ we have as a consequence of
Theorem \ref{tabonadic} that $(x,S,a)$ is part of a 
tabulation of $ax^*\f X\ra A$. It follows that we have a natural
isomorphism $(x,S,a)\ra FC(x,S,a)$, which completes the demonstration
that $C_{X,A}$ and $F_{X,A}$ are inverse equivalences.
\frp

\section{Direct sums in bicategories of spans}

In the previous section we gave a characterization of those (cartesian)
bicategories of the form $\spn\E$ for a category $\E$ with finite
limits. In this final section we give a refinement, showing that $\spn\E$
has direct sums if and only if the original category $\E$ is lextensive \cite{ext}.

Direct sums are of course understood in the bicategorical sense. A {\em zero object}
in a bicategory is an object which is both initial and terminal. In a bicategory with
finite products and finite coproducts in which the initial object is also terminal
there is a canonical induced arrow $X+Y\to X\times Y$, and we say that the 
bicategory has {\em direct sums} when this map is an equivalence. 

\begin{remark}
Just as in the case of ordinary categories, the existence of direct sums gives rise
to a calculus of matrices. A morphism $X_1+\ldots+X_m\to Y_1+\ldots+Y_n$
can be represented by an $m\times n$ matrix of morphisms between the summands,
and composition can be represented by matrix multiplication.
 \end{remark}

\thm
Let $\E$ be a category with finite limits, and $\bB=\spn\E$. Then the 
following are equivalent:
\begin{enumerate}[$i)$]
\item $\bB$ has direct sums;
\item $\bB$ has finite coproducts;
\item $\bB$ has finite products;
\item $\E$ is lextensive.
\end{enumerate}
\eth

\prf

$[i)\Longrightarrow$ $ii)]$ is trivial.

$[ii)\Longleftrightarrow$ $iii)]$ follows from the fact that $\bB\op$ is 
biequivalent to $\bB$.

$[ii)\Longrightarrow$ $iv)]$  
Suppose that $\bB$
has finite coproducts, and write $0$ for the initial object and $+$ for the coproducts. 

For every object $X$ there is a unique span $0\la D\ra X$. By uniqueness, any
map into $D$ must be invertible, and any two such with the same domain must
be equal. Thus when we compose the span with its opposite, as in 
$0\la D\ra X\la D\ra 0$, the resulting span is just $0\la D\ra 0$. Now by the 
universal property of $0$ once again, this must just be $0\la 0\ra 0$, and so 
$D\cong 0$, and our unique span $0\to X$ is a map.

Clearly coproducts of maps are coproducts, and so the coproduct injections
$X+0\to X+Y$ and $0+Y\to X+Y$ are also maps. Thus the coproducts in $\bB$
will restrict to $\E$ provided that the codiagonal $u\f X+X\la E \ra X\f v$ is a map for
all objects $X$. Now the fact that the codiagonal composed with the first injection
$i:X\to X+X$ is the identity tells us that we have a diagram as on the left below
$$
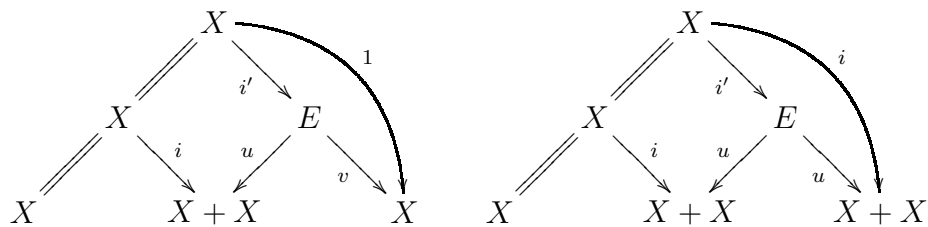
\xymatrix @!R=1pc @!C=1pc {
&& X \ar[dr]_{i'} \ar@{=}[dl] \ar@/^2pc/[ddrr]^{1} && & 
&& X \ar[dr]_{i'} \ar@{=}[dl] \ar@/^2pc/[ddrr]^{i} \\
& X \ar[dr]^{i} \ar@{=}[dl] && E \ar[dl]_{u} \ar[dr]_{v} & & 
& X \ar[dr]^{i} \ar@{=}[dl] && E \ar[dl]_{u} \ar[dr]_{u} \\
X && X+X && X & X && X+X && X+X }$$
in which the square is a pullback; but then the diagram on the right 
shows that the composite of $u\f X+X\la E\ra X+X\f u$ with the injection
$i:X\to X+X$ is just $i$. Similarly its composite with the other injection
$j:X\to X+X$ is $j$, and so $u\f X+X\la E\ra X+X\f u$ is the identity. 
This proves that the codiagonal is indeed a map, and so that $\E$ has finite
coproducts; we have already assumed that it has finite limits. To see that $\E$ 
is lextensive observe that we have equivalences 
$$\E/(X+Y)\simeq \bB(X+Y,1) \simeq \bB(X,1)\times\bB(Y,1)\simeq \E/X\times \E/Y.$$

$[iv)\Longrightarrow$ $i)]$  
Suppose that $\E$ is lextensive. Then in particular, it is distributive, so that
$(X+Y)\x Z\cong X\x Z+X\x Y$, and 
we have 
\begin{align*}
\bB(X+Y,Z) &\simeq \E/\bigl((X+Y)\x Z\bigr) \simeq \E/(X\x Z+Y\x Z) \\
                  &\simeq \E/(X\x Z)\times \E/(Y\x Z) \simeq \bB(X,Z)\times\bB(Y,Z)
\end{align*}
which shows that $X+Y$ is the coproduct in $\bB$; but a similar argument shows
that it is also the product.
\frp

\rmk
The implication $iv)\Rightarrow i)$ was proved in \cite[Section~3]{SP07}.
\eth

\rmk
The equivalence $ii)\Leftrightarrow iv)$ can be seen as a  special case of a 
more general result \cite{HS} characterizing colimits in $\E$ which are also (bicategorical) colimits in $\spn\E$.

\rmk
There is a corresponding result involving partial maps in lextensive categories,
although the situation there is more complicated as one does not have direct
sums but only a weakened relationship between products and coproducts, and 
a similarly weakened calculus of matrices. See \cite[Section~2]{restiii}.
\eth

There is also a nullary version of the theorem. We simply recall that an initial object in a
category $\E$ is said to be {\em strict}, if any morphism into it is invertible, and 
then leave the proof to the reader. Once again the equivalence $ii)\Leftrightarrow iv)$
is a special case of \cite{HS}.

\thm
Let $\E$ be a category with finite limits, and $\bB=\spn\E$. Then the 
following are equivalent:
\begin{enumerate}[$i)$]
\item $\bB$ has a zero object;
\item $\bB$ has an initial object;
\item $\bB$ has a terminal object;
\item $\E$ has a strict initial object.
\end{enumerate}
\eth

\end{document}